\documentclass[final]{siamltex}

\usepackage{amsfonts}
\usepackage{amsmath}
\usepackage[english]{babel}
\usepackage{subfigure}
\usepackage{enumerate}
\usepackage{graphicx}

\pdfoutput=1

\newcommand{\CC}{\mathbb{C}}

\newcommand{\RR}{\mathbb{R}}

\newcommand{\ints}{\int\limits}

\newcommand{\lam}{\lambda}

\newcommand{\vt}{\vartheta}
\newcommand{\OO}{\mathcal{O}}
\newcommand{\de}{\partial}

\newcommand{\zi}{\zeta}
\newcommand{\somma}{\sum\limits_{j\in\{ s,a \}}}
\newcommand{\ja}{j\in\{ s,a \}}
\newcommand{\Cam}{\mathcal{C}}
\newcommand{\lsomma}{\sum\limits_{l=0}^{M}}
\newcommand{\vp}{\varphi}

\DeclareMathOperator{\sign}{\mathrm{sign}}
\DeclareMathOperator{\RE}{\mathrm{Re}}
\DeclareMathOperator{\IM}{\mathrm{Im}}
\DeclareMathOperator{\arcsinh}{\mathrm{arcsinh}}

\newtheorem{remark}[theorem]{Remark}

\title{A radiation condition for uniqueness in a wave propagation problem for 2-D open waveguides}

\author{Giulio Ciraolo \thanks{Dipartimento di Matematica e
Applicazioni per l'Architettura, Universit\`a di Firenze, Piazza
Ghiberti 27, 50122 Firenze, Italy, ({\tt ciraolo@math.unifi.it}).}
\and Rolando Magnanini \thanks{Dipartimento di matematica ``U.
Dini'', Universit\`a di Firenze, Viale Morgagni 67/A, 50134 Firenze,
Italy, ({\tt magnanin@math.unifi.it}).}}

\begin{document}

\maketitle

\begin{abstract}
We study the uniqueness of solutions of Helmholtz equation for a
problem that concerns wave propagation in waveguides. The classical
radiation condition does not apply to our problem because the
inhomogeneity of the index of refraction extends to infinity in one
direction. Also, because of the presence of a waveguide, some waves
propagate in one direction with different propagation constants and
without decaying in amplitude.

Our main result provides an explicit condition for uniqueness which
takes into account the physically significant components,
corresponding to guided and non-guided waves; this condition reduces
to the classical Sommerfeld-Rellich condition in the relevant cases.

Finally, we also show that our condition is satisfied by a solution,
already present in literature, of the problem under consideration.
\end{abstract}



\pagestyle{myheadings} \thispagestyle{plain} \markboth{G. CIRAOLO
AND R. MAGNANINI}{UNIQUENESS FOR 2-D WAVEGUIDES}

\section{The problem of uniqueness for the Helmholtz equation} \label{section introd
rad cond} Let $\Sigma \subset \RR^N$ be a (possibly empty) bounded
closed surface. It is well known that the Dirichlet problem
\begin{equation} \label{Dirichlet pb}
\begin{cases}
\Delta u + k^2 u = f & \textmd{outside } \Sigma ,\\
u = U & \textmd{on } \Sigma ,
\end{cases}
\end{equation}
has not an unique solution. If $k=0$ (Poisson's equation), in order
to obtain the uniqueness, it is required that the solution vanishes
at infinity. If $k \neq 0$, that is not sufficient anymore. In fact,
there are two different solutions of \eqref{Dirichlet pb} which
vanish at infinity, representing the outward and inward radiation.
Hence, an additional (or different) condition at infinity is needed.

The first condition we can add is that
\begin{equation} \label{sommerfeld}
  \lim\limits_{R\to \infty} R^{\frac{N-1}{2}} \left( \frac{\de u}{\de R} - ik u \right)=0,
\end{equation}
uniformly; this is the so-called \emph{Sommerfeld's radiation
condition}. Here, $\frac{\partial u}{\partial R}$ denotes the radial
derivative of $u$. The physical meaning of this condition is that
there are no sources of energy  at infinity. Moreover, it assures
that, far from the surface $\Sigma$, $u$ behaves as a wave generated
by a point source.

Stated as in \eqref{sommerfeld} and together with the assumption
that $u$ vanishes at infinity, this condition is due to Sommerfeld,
see \cite{So1} and \cite{So2} (see also \cite{Mag1} and
\cite{Mag2}). The vanishing assumption on $u$ was dropped by Rellich
(see \cite{Rel}), who also proved that \eqref{sommerfeld} can be
replaced by the weaker condition
\begin{equation}\label{rellich surface}
\lim\limits_{R\to \infty} \ints_{\partial B_R} \left| \frac{\de
u}{\de R} - ik u \right|^2 d \sigma = 0 \ ,
\end{equation}
where $B_R$ is the ball centered at the origin with radius $R$. In
the same paper, Rellich also proved that a radiation condition can
also be given in the form
\begin{equation}\label{rellich}
\ints_{\RR^N} \left| \frac{\de u}{\de R} - ik u \right|^2 dx <
+\infty.
\end{equation}
Condition \eqref{rellich} can be considered the starting point for
our work, as we are going to explain shortly. Before describing our
results, we cite some generalizations of the work of Rellich.

When $n$ is a function which is identically $1$ outside a compact
set, \eqref{rellich surface} still guarantees the uniqueness of a
solution of
\begin{equation*}
\Delta u + k^2 n(x)^2 u = f, \quad x\in \RR^N;
\end{equation*}
see \cite{Mi1} and \cite{Sc} and references therein.

Several authors (see \cite{Rel} \cite{Mi2} \cite{RS} \cite{Zh}
\cite{PV} \cite{Ei}) studied the case in which $n$ is not constant
at infinity, but has an angular dependency, say $n(x) \to
n_\infty(x/ |x|)$ as $|x|\to\infty$, and it approaches to the limit
with a certain behaviour. Among these papers, we want to mention the
results in \cite{Zh} and \cite{PV}, where the authors proved the
uniqueness of solutions of the Helmholtz equation by means of the
limiting absorbtion method and by introducing the radiation
condition:
\begin{equation*}
\lim_{R\to +\infty} \frac{1}{R} \ints_{B_R} \Big{|} \frac{\partial
u}{\partial R} - i k n_{\infty} u \Big{|}^2 dx = 0.
\end{equation*}
Here, the assumptions on $n$ are so that the energy cannot be
trapped along any direction, but it radiates toward infinity. That
is in accordance with \cite{PS}, where the authors point out that
the Sommerfeld radiation condition, since it involves the dimension,
is inappropriate for problems admitting a lower dimensional solution
(a plane wave).

The present paper is motivated by the study of wave propagation in
optical waveguides. In particular, we shall study the uniqueness of
solutions of the two-dimensional Helmholtz equation
\begin{equation}\label{helm}
\Delta u + k^2 n(x)^2 u = f, \quad (x,z) \in \RR^2,
\end{equation}
where $n$ is of the form
\begin{equation}\label{n}
  n:= \begin{cases} n_{co}(x), & |x|\leq h,\\
  n_{cl}, & |x| > h; \end{cases}
\end{equation}
here $n_{co}$ is a bounded function and $n_{cl}$ is a constant;
thus, (\ref{n}) models the index of refraction of a rectilinear open
waveguide of width $2h$ (subscripts $co$ and $cl$ refer to the
{\it core} and {\it cladding} of the waveguide) (see \cite{SL}).

We observe that functions $n$ of type (\ref{n}) are not considered
in the works cited before. In fact, the most important feature of
optical waveguides is the presence of waves confined inside the
waveguide (\emph{guided modes}) which are oscillatory and never
decaying along the axis of propagation ($z$-axis). It is easy to
show that a pure guided mode supported by the Helmholtz equation
does not satisfy the radiation conditions above retrieved (as
already pointed out in \cite{PS}). Functions $n$ similar to
\eqref{n} were considered by J\"{a}ger and Sait\={o} in
\cite{JS1}-\cite{JS2}; however, their assumptions on $n$ do not
admit the occurrence of guided modes.

As far as we know, the only works dealing with uniqueness in an
optical waveguide setting have appeared in the Russian literature
(see \cite{Rei} \cite{No} \cite{NS} \cite{KNH} and references
therein). However, the {\it Reichardt condition} studied therein
only deals with guided modes and does not apply to the total field.

The main result of this paper is Theorem \ref{teo unicita}, where we
present a new radiation condition that guarantees the uniqueness of
a solution of \eqref{helm} with $n$ given by \eqref{n}. We observe
that, if we suppose that no guided mode is present (this is possible
by choosing special parameters in the function $n$), our radiation
condition reduces to \eqref{rellich}. In this setting, our results
provide a different proof of special cases studied in \cite{Rel} and
\cite{JS1}-\cite{JS2}.

The key ingredients of our proof are essentially four: (i) if
\eqref{helm} possesses two solutions satisfying our radiation
condition, then their difference $w$ must belong to the Sobolev
space $H^2(\RR^2)$; (ii) as a consequence of (i), the Fourier
transform of $w$ in the $z$-direction (parallel to the fiber's axis)
is square integrable for almost all $x\in\RR$ and satisfies an
ordinary differential equation in $x$; (iii) the only square
integrable solution of such an equation is identically zero; (iv)
the proof is then completed by using an appropriate transform theory
in the $x$-direction and repeating the arguments in (ii) and (iii).
This scheme will be carried out in \S \ref{section 2}.

In \cite{MS} the authors derived a solution\footnote{We will refer
to such a solution as the \emph{spectrum-based solution}.} for the
problem \eqref{helm}-\eqref{n} in terms of a Green's function.
Section \ref{section 3} is devoted to prove that such a solution
satisfies our radiation condition. This will be done in three steps:
in \S \ref{subsection 3.1} we derive a representation of the
solution as a contour integral; in \S \ref{subsection 3.2} we prove
uniform estimates for the non-guided part of the
\emph{spectrum-based solution}; in \S \ref{subsection 3.3} we carry
out the proof by testing the radiation condition on the guided part
and using the asymptotic estimates obtained in \S \ref{subsection
3.2}.

We wish to observe that the results in the present paper can be
easily adapted to prove the uniqueness of a solution for the
Pekeris waveguide problem (see \cite{Wi} and \cite{dS}).

\section{A new Rellich-type condition and a uniqueness theorem}
\label{section 2} In this section we shall state a radiation
condition, that generalizes \eqref{rellich}, and prove our uniqueness
result.

\subsection{Preliminaries}\label{section prelimin}
We recall the relevant results of \cite{MS} which will be useful in the
rest of the paper.

In \cite{MS} a Green's
function $G$ for \eqref{helm} is constructed:
a solution of \eqref{helm} is
\begin{equation} \label{u}
u(x,z) = \ints_{\RR^2} G(x,z;\xi,\zi) f (\xi,\zi) d\xi d\zi,
\end{equation}
where
\begin{equation}\label{Green ch2}
  G(x,z;\xi,\zi) = \sum_{j\in\{s,a\}} \ints_{0}^{+\infty} \frac{e^{i|z-\zi| \sqrt{k^2n_*^2 - \lam}}}{2i
  \sqrt{k^2n_*^2 -\lam}} v_j(x,\lam) v_j(\xi,\lam) d\rho_j(\lam).
\end{equation}
Here
\begin{equation*} 
n_*=\max_{\RR}{n},
\end{equation*}
\begin{equation*}
  \langle d\rho_{j}, \eta \rangle = \sum_{m=1}^{M_j} r_m^j \eta (\lam_m^j) + \frac{1}{2\pi} \ints_{d^2}^{+\infty}
  \frac{\sqrt{\lam -d^2}}{(\lam - d^2) \phi_j(h,\lam)^2 + \phi_j'(h,\lam)^2} \eta(\lam)
  d\lam,
\end{equation*}
for all $\eta\in C_0^\infty(\RR)$ ($\langle\:,\:\rangle$ is the
usual dual product),
\begin{equation}\label{rmj}
  r_m^j = \left[ \ints_{-\infty}^{+\infty} v_j(x,\lam_m^j)^2 dx \right]^{-1} =
  \frac{\sqrt{d^2-\lam_m^j}}{ \sqrt{d^2-\lam_m^j} \ints_{-h}^h \phi_j(x,\lam_m^j)^2 dx +
  \phi_j(h,\lam_m^j)^2}.
\end{equation}
Also, $v_j(x,\lam)$ are linearly independent solutions of
\begin{equation}\label{eq v ch2}
  v'' + [\lam - q(x)] v = 0, \quad \textmd{in } \RR,
\end{equation}
where $q(x)=k^2 [n_*^2-n(x)^2]$, and have the following form:
\begin{equation}\label{vj ch2}
v_j(x,\lam)=
  \begin{cases}
    \phi_j(h,\lam) \cos Q (x-h) + \frac{\phi_j'(h,\lam)}{Q} \sin Q(x-h), & \text{if } x > h, \\
    \phi_j(x,\lam), & \text{if } |x| \leq h, \\
    \phi_j(-h,\lam) \cos Q (x+h) + \frac{\phi_j'(-h,\lam)}{Q} \sin Q(x+h), & \text{if } x < -h, \\
  \end{cases}
\end{equation}
for $j=s,a$, with $Q=\sqrt{\lam-d^2}$, $d^2=k^2(n_*^2-n_{cl}^2)$;
the $\phi_j$'s are solutions of \eqref{eq v ch2} in the interval
$(-h,h)$ and satisfy the initial conditions:
\begin{equation*}
  \begin{array}{cc}
    \phi_s(0,\lam)=1, & \phi_s'(0,\lam)=0, \\
    \phi_a(0,\lam)=0, & \phi_a'(0,\lam)=\sqrt{\lam}.
  \end{array}
\end{equation*}
(The indices $j=s,a$ correspond to symmetric and antisymmetric
solutions, respectively.)

We notice that \eqref{Green ch2} can be split up into two summands,
\begin{equation*}
  G=G^g + G^{rad},
\end{equation*}
where
\begin{equation}\label{G^g ch2}
 G^g(x,z;\xi,\zi) = \sum_{j\in \{ s,a \}} \sum_{m=1}^{M_j} \frac{e^{i|z-\zi|
 \sqrt{k^2 n_*^2 - \lam_m^j}}}{2i \sqrt{k^2n_*^2 -\lam_m^j}}
 v_j(x,\lam_m^j) v_j(\xi,\lam_m^j) r_m^j,
\end{equation}
and
\begin{equation}\label{G^rad ch2}
  G^{rad} (x,z;\xi,\zi) = \frac{1}{2\pi}  \sum_{j\in\{s,a\}} \ints_{d^2}^{+\infty}
  \frac{e^{i|z-\zi| \sqrt{k^2n_*^2 - \lam}}}{2i \sqrt{k^2n_*^2 -\lam}}
  v_j(x,\lam) v_j(\xi,\lam) \frac{\sigma_j(\lam)}{\sqrt{\lam-d^2}} d\lam,
\end{equation}
with
\begin{equation}\label{sigma ch2}
\sigma_j(\lam) = \frac{\lam -d^2}{(\lam - d^2) \phi_j(h,\lam)^2 +
\phi_j'(h,\lam)^2},\quad j=s,a;
\end{equation}

$G^g$ represents the guided part of the Green's function, which
involves the guided modes, i.e. the modes propagating mainly inside
the waveguide; $G^{rad}$ is the part of the Green's function
corresponding to the non-guided energy, i.e. the energy radiated
outside or vanishing along the waveguide, which we denote by
\begin{equation}\label{u rad originale}
u^{rad}= \ints_{\RR^2} G^{rad}(x,z;\xi,\zi) f(\xi,\zi).
\end{equation}

It exists a finite number of guided modes, which corresponds to the
finite number of roots of the equations
\begin{equation*}
\sqrt{d^2-\lam} \; \phi_j(h,\lam) + \phi_j'(h,\lam) = 0, \quad
j\in\{ s,a \},
\end{equation*}
laying in the interval $(0,d^2)$. We shall denote by $\lam_m^j$,
$m=1,\ldots, M_j$, $j=s,a$, these roots. Each $v_j(x,\lam_m^j)$
decays exponentially for $|x| > h$ as it is clear from the formula:
\begin{equation}\label{vj guidate}
v_j(x,\lam_m^j) = \begin{cases} \phi_j(h,\lam_m^j) e^{- \sqrt{d^2 -\lam_m^j} (x-h)}, & x > h, \\
\phi_j(x,\lam_m^j), & |x|\leq h,\\
\phi_j(-h,\lam_m^j) e^{\sqrt{d^2 -\lam_m^j} (x+h)}, & x < - h.
\end{cases}
\end{equation}
We notice that $G^g$ is bounded and oscillatory in the $z$
direction, because $\sqrt{k^2 n_*^2 - \lam_m^j}$ is real for every
$m=1,\ldots, M_j$, $j=s,a$.

\vspace{1em}

\begin{remark}\label{remark meromorphic}
The functions $\sigma_j(\lam)$, $\ja$, given by \eqref{sigma ch2},
are meromorphic functions of $\lam\in\CC$, real-valued for
$\lambda\in\RR$ and with poles that are real and simple (see
\cite{CL},\cite{Ti}), which corresponds to the values $\lam_m^j$,
$m=1,\ldots, M_j$, $j=s,a$.
\end{remark}

\vspace{1em}

To simplify notations, we shall denote by $\gamma_l$,
$l=1,\ldots,M,\ M=M_s+M_a$, the values $\lam_m^j,\ m=1,\ldots,M_j,\
j=s,a$, ordered according to the natural ordering on the real line,
and by $\gamma_*$ their maximum. With these premises, we set
\begin{equation}\label{e autof}
e(x,\gamma_l)=\frac{v_j(x,\gamma_l)}{\| v_j(\cdot,\gamma_l)\|_2}.
\end{equation}
From \eqref{G^g ch2} and \eqref{rmj},
it is clear that the guided part $G^{g}$ can be written as
\begin{equation*}
G^g(x,z;\xi,\zi)= \sum_{l=1}^M G_l^g(x,z;\xi,\zi),
\end{equation*}
where
\begin{equation}\label{G l guid}
G_l^g(x,z;\xi,\zi) = \frac{e^{i\beta_l |z-\zi|}}{2i \beta_l}
e(x,\gamma_l) e(\xi,\gamma_l),
\end{equation}
with
\begin{equation}\label{beta_l}
\beta_l= \sqrt{k^2 n_*^2 - \gamma_l}.
\end{equation}

\vspace{1em}

Let $s\in\RR;$ we will denote by $L^{2,s}(\RR^2)$ the weighted
Lebesgue space consisting of all the complex-valued measurable
functions $u$ such that $(1+x^2+z^2)^s |u(x,z)|^2$ is summable in
$\RR^2,$ equipped with the natural norm defined by
\begin{equation*}
|u|_{2,s}^2 = \ints_{\RR^2} |u(x,z)|^2 (1+x^2+z^2)^s dx dz;
\end{equation*}
$L^{2,s}(\RR^2)$ is commonly used when dealing with solutions of
Helmholtz equation (see \cite{Ag} and \cite{Le}). In \cite{CM} we
proved that the spectrum-based solution (\ref{u})-(\ref{Green ch2})
derived in \cite{MS} belongs to $L^{2,s}(\RR^2)$, for $s<-1$, if
$f\in L^{2,-s}(\RR^2)$.

\vspace{1em}

The following lemma will be useful in the next subsection.

\vspace{1em}

\begin{lemma}\label{lemma crescita}
For $s<-1$, let $w \in L^{2,s}(\RR^2)$ satisfy
\begin{equation}
\label{helm omogenea}
\Delta u + k^2 n(x)^2 u = 0
\end{equation}
in $\RR^2,$ where $n$ is given by (\ref{n}).
Then
\begin{equation}\label{H1}
\lim\limits_{|x| \to +\infty} u  (x,z) e^{-|x| \sqrt{d^2 -
\gamma_*}}=\lim\limits_{|x| \to +\infty} u_x (x,z)
e^{-|x|\sqrt{d^2 - \gamma_*}} = 0,
\end{equation}
where $\gamma_* = \max\limits_{1\le l\le M}{\gamma_l}$.
\end{lemma}

\begin{proof}
Since $u$ is a solution of \eqref{helm omogenea}, from Lemmas A.1
and A.3 in \cite{CM}, we infer that both $(1+x^2+z^2)^s |\nabla u(x,z)|^2$
and $(1+x^2+z^2)^s |\nabla^2 u(x,z)|^2$ are summable in $\RR^2$. Thus, it
easily follows that the function
$$
\Psi(x,z) = (1+x^2+z^2)^{s/2} u(x,z)
$$
belongs to the Sobolev space $W^{2,2}(\RR^2)$. The Sobolev Imbedding Theorem (see
Theorem 4.12 in \cite{AF}) implies that $\Psi \in L^\infty (\RR^2)$
and hence the first limit in \eqref{H1} follows at once.

A straightforward computation shows that $\Psi$ satisfies the following equation
$$
\Delta\Psi+b\cdot\nabla\Psi+c\ \Psi=0
$$
in $\RR^2,$ where
$$
b(x,z)=-\frac{2s (x,z)}{1+x^2+z^2}, \quad c(x,z)=k^2 n(x)^2-s\ \frac{2-s(x^2+z^2)}{(1+x^2+z^2)^2}.
$$
Since $\Psi\in W^{2,2}(\RR^2),$ by Theorem 8.10 in \cite{GT}, we have that $\Psi\in W^{3,2}(H_+)$
where $ H_+ = \{ (x,z) \in \RR^2 : x \geq h \}$. Again, by the Sobolev
Imbedding Theorem, $|\nabla \Psi|$ is bounded in $H_+$ and hence
the second limit in \eqref{H1} holds as $x\to +\infty$.
The same limit as $x\to-\infty$ holds by a similar argument.
\end{proof}

\subsection{The radiation condition and uniqueness theorem}
We consider a solution $u$ of \eqref{helm} and define
\begin{equation}\label{u_l}
    u_l(x,z)= e(x,\gamma_l) U(z,\gamma_l), \quad l=1,\ldots,M,
\end{equation}
with $e(x,\gamma_l)$ given by \eqref{e autof} and where
\begin{equation} \label{U def}
U(z,\gamma_l)= \ints_{-\infty}^\infty u(\xi,z) e(\xi,\gamma_l) d\xi,
\quad l=1,\ldots, M.
\end{equation}
The remainder part of $u$ is
\begin{equation}\label{u_0}
    u_0 (x,z)= u(x,z) - \sum\limits_{l=1}^{M}  u_l(x,z).
\end{equation}
%

We introduce a one-parameter family of sets $\Omega_R$, $R>0$, such
that $\de \Omega_R$ are level sets of a convex and coercive function
$d(x,z)$, i.e. $\Omega_R=\{(x,z)\in \RR^2 : \ d(x,z)\leq R \}$.

With these notations, we state our {\it radiation condition} for a
solution $u$ of \eqref{helm}:
\begin{equation}\label{rad cond}
    \lsomma \ints_0^\infty \ints_{\partial \Omega_R}
    \Big{|} \frac{\partial u_l}{\partial \nu} - i\beta_l u_l
    \Big{|}^2 d\ell \, d R < + \infty,
\end{equation}
with $\beta_0 = k n_{cl}$ and $\beta_l$,$l=1,\ldots,M$, given by
\eqref{beta_l}.

Notice that, when $n\equiv 1$, we can choose $\Omega_R=B_R$ and
\eqref{rad cond} reduces to \eqref{rellich}, since in such a case
the guided components are not present.

\vspace{1em}

The main result of this paper follows.
\vspace{1em}
\begin{theorem}\label{teo unicita}
There is at most one solution  of
\eqref{helm} that satisfies \eqref{rad cond}
and belongs to $u\in L^{2,s}(\RR^2),$ $s<-1.$
\end{theorem}
\vspace{1em}
\begin{remark}\label{remark 1}
As it will be clear, it is not necessary to specify further the sets
$\Omega_R$ in \eqref{rad cond} to get uniqueness of a solution of
\eqref{helm}. That means that Theorem \ref{teo unicita} holds for
{\it any} choice of one-parameter family of sets $\Omega_R$
satisfying the above mentioned assumptions.

Of course, a solution of \eqref{helm} satisfying \eqref{rad cond}
may not exist for an arbitrary choice of the sets $\Omega_R.$ In \S
3, we shall choose a special family of sets $\Omega_R$ and prove
that, with this choice, the solution of \eqref{helm} given by
\eqref{u}-\eqref{Green ch2} satisfies \eqref{rad cond}.

We also notice that it is not necessary to choose the same sets
$\Omega_R$ in each addendum in \eqref{rad cond}; Theorem \ref{teo unicita} still holds
if we replace \eqref{rad cond} by the more general radiation condition
\begin{equation}\label{rad cond 2}
\lsomma \ints_0^\infty \ints_{\partial \Omega_R^{(l)}}
    \Big{|} \frac{\partial u_l}{\partial \nu} - i\beta_l u_l
    \Big{|}^2 d\ell \, d R < + \infty,
\end{equation}
where $\Omega_R^{(l)},$ $l=0,1,\dots, M,$  are one-parameter
families satisfying the above mentioned assumptions.
\end{remark}

\vspace{1em}

\noindent Theorem \ref{teo unicita} is based on Lemma \ref{lemma
parte immaginaria 0} and Theorem \ref{teo wm L2} below.

\vspace{1em}

\begin{lemma} \label{lemma parte immaginaria 0}
Let $\beta \in \RR$ and $u$ be a weak solution of \eqref{helm omogenea}
Then
\begin{equation} \label{eq parte immag 0}
\ints_{\partial \Omega} \Big{|} \frac{\partial w}{\partial \nu} -
i\beta w \Big{|}^2 d\ell = \ints_{\partial \Omega} \left( \Big{|}
\frac{\partial w}{\partial \nu} \Big{|}^2 + \beta^2 |w|^2 \right)
d\ell,
\end{equation}
for every $\Omega \subset \RR^2$ bounded and sufficiently smooth.
\end{lemma}

\begin{proof}
Since $u$ is a weak solution of \eqref{helm omogenea}, by Theorem
8.8 in \cite{GT}, we obtain the necessary regularity to infer that
\begin{eqnarray*}
&&\ints_{\partial \Omega} \bar{u}\ \frac{\de u}{\de \nu}\ d\ell =
\ints_\Omega \mbox{\rm div}(\bar{u}\nabla u)\ dx dz=
\ints_\Omega \left\{|\nabla u|^2+\bar{u} \Delta u\right\} dx dz=\\
&&\ints_\Omega \left\{|\nabla u|^2-k^2 n(x)^2 |u|^2\right\}\ dx dz.
\end{eqnarray*}

Therefore
\begin{equation*}
\IM \ints_{\de \Omega} \bar{w}\ \frac{\de w}{\de \nu}\ d\ell = 0,
\end{equation*}
which easily implies \eqref{eq parte immag 0}.
\end{proof}

\vspace{1em}

\begin{theorem}\label{teo wm L2}
Let $u\in L^{2,s} (\RR^2)$ be a weak solution of \eqref{helm omogenea}
satisfying \eqref{rad cond}. Then
\begin{equation}\label{w_m w mu in L2}
\lsomma \ \ints_0^{+\infty} dR \ints_{\partial \Omega_R} \left[
\Big{|} \frac{\partial u_l}{\partial \nu} \Big{|}^2 + \beta_l^2
|u_l|^2 \right] d\ell < + \infty,
\end{equation}
and, in particular,
\begin{equation} \label{wm norma finita}
\ints_{\RR^2} |u_l|^2 dx dz < +\infty,
\end{equation}
for every $l=0,1,\ldots,M$.
\end{theorem}

\begin{proof}
By Lemma \ref{lemma parte immaginaria 0}, it is enough to prove that
each $u_l$, $l=0,1,\ldots,M$, satisfies \eqref{helm omogenea}. Then,
\eqref{w_m w mu in L2} and \eqref{wm norma finita} will follow from
\eqref{eq parte immag 0} and \eqref{rad cond}.

Suppose $l\geq 1$. Since
\begin{equation} \label{helm weak omogenea}
-\ints_{\RR^2} \nabla u \cdot \nabla \varphi \; dx dz + k^2
\ints_{\RR^2} n(x)^2 u \varphi \; dx dz =0
\end{equation}
for every $\varphi \in H_0^1(\RR^2)$, we choose
$\vp (x,z) = e(x,\gamma_l) \eta(z)$ with $\eta \in C_0^1 (\RR)$, and obtain:
\begin{equation*}
\begin{split}
-\ints_{\RR^2}[ u_x(x,z) & e'(x,\gamma_l) \eta(z)  +u_z(x,z) e(x,\gamma_l) \eta'(z)]\ dx dz +\\
&k^2 \ints_{\RR^2} n(x)^2 u(x,z) e(x,\gamma_l) \eta(z) dx dz = 0;
\end{split}
\end{equation*}
an integration by parts and Lemma \ref{lemma crescita} then give
\begin{equation*}
\begin{split}
\ints_{\RR^2} u(x,z) &  e''(x,\gamma_l) \eta(z) dx dz - \ints_{\RR^2} u_z(x,z) e(x,\gamma_l) \eta'(z) dx dz +\\
&k^2 \ints_{\RR^2} n(x)^2 u(x,z) e(x,\gamma_l) \eta(z) dx dz = 0.
\end{split}
\end{equation*}
Since $e(x,\gamma_l)$ satisfies \eqref{eq v ch2}, we obtain
\begin{equation*}
- \ints_{\RR^2} u_z(x,z) e(x,\gamma_l) \eta'(z) dx dz + (k^2 n_*^2 - \gamma_l)
\ints_{\RR^2} u(x,z) e(x,\gamma_l) \eta(z) dx dz = 0,
\end{equation*}
and thus, from \eqref{U def},
\begin{equation*}
-\ints_{\RR} U_z (z,\gamma_l) \eta'(z) dz + (k^2 n_*^2 - \gamma_l)
\ints_{\RR} U(z,\gamma_l) \eta(z) dz = 0,
\end{equation*}
for every $\eta \in C_0^1(\RR)$.

Together with \eqref{eq v ch2}, this formula implies that each
$u_l(x,z)$, $l=1,\ldots,M$, given by \eqref{u_l},
is a weak solution of \eqref{helm
omogenea}. In fact, for $\vp (x,z) = \psi(x) \eta(z)$ with
$\psi,\eta \in C_0^1 (\RR),$ integration by parts gives
\begin{equation*}
\begin{split}
-\ints_{\RR^2} &\nabla u_l(x,z)  \cdot \nabla \varphi(x,z)\ dx dz +  k^2
\ints_{\RR^2} n(x)^2 u_l(x,z) \varphi(x,z)\ dx dz= \\
& \left(\ints_{\RR}\{e''(x,\gamma_l)+ [\gamma_l - q(x)]e(x,\gamma_l)\} \psi(x)\ dx\right) \left(\ints_\RR U(z,\gamma_l) \eta(z) dz\right) +\\
&\left(\ints_{\RR} e(x,\gamma_l) \psi(x) dx\right)\left(\ints_{\RR} [-U_z(z,\gamma_l)\eta'(z) + (k^2 n_*^2 - \gamma_l)
U(z,\gamma_l) ]\eta(z) dz\right)=0  \: ;
\end{split}
\end{equation*}
the same conclusion holds for any $\vp\in C^1_0(\RR^2),$ by a density argument.

Since $u$ and $u_l$, $l=1,\ldots,M$, now satisfy \eqref{helm omogenea},
the same holds for $u_0$. Thus, as already mentioned, we can apply
Lemma \ref{lemma parte immaginaria 0} to each $u_l$,
$l=0,1,\ldots,M$, and obtain
\begin{equation*}
\lsomma \ \ints_{\partial \Omega_R} \Big{|} \frac{\partial
u_l}{\partial \nu} - i\beta_l u_l \Big{|}^2 d\ell = \lsomma \
\ints_{\partial \Omega_R} \left( \Big{|} \frac{\partial
u_l}{\partial \nu} \Big{|}^2 + \beta_l^2 |u_l|^2 \right) d\ell,
\end{equation*}
for every $R>0$, and then, since $u$ satisfies \eqref{rad cond}, we
get \eqref{w_m w mu in L2} and \eqref{wm norma finita}.
\end{proof}


\subsection{Proof of Theorem \ref{teo unicita}}
Let $u_1$ and $u_2$ be two solutions; $u=u_1 - u_2$
satisfies \eqref{helm omogenea} and \eqref{rad cond}.

From Theorem \ref{teo wm L2} we have that $u\in L^2(\RR^2)$ and, by
using Lemmas A.1 and A.3 in \cite{CM}, we get $u\in H^2(\RR^2)$.
Therefore, $u(x,\cdot) \in L^2(\RR)$ for almost every $x\in\RR$, and
the same holds for $u_x(x,\cdot)$ and $u_{xx}(x,\cdot)$. Hence, we
can transform \eqref{helm omogenea} by using the Fourier transform
in the $z$-coordinate,
\begin{equation*}
\hat{u}(x,t) = \ints_{-\infty}^{+ \infty} u(x,z) e^{-izt} dz, \quad
\textmd{for a.e. } x\in\RR,
\end{equation*}
and obtain:
\begin{equation}\label{helm fourier in z}
\hat{u}_{xx} (x,t) + [k^2 n(x)^2 - t^2] \hat{u} (x,t) = 0, \quad
\textmd{a.e. } x\in\RR.
\end{equation}

From Fubini-Tonelli's theorem, the integrals
\begin{equation*}
\ints_{\RR^2} |\hat{u}(x,t)|^2 dx dt,\ \ \ints_{-\infty}^{+\infty}
dt \ints_{-\infty}^{+\infty}|\hat{u}(x,t)|^2 dx\ \ \textmd{and }\
\ints_{-\infty}^{+\infty} dx \ints_{-\infty}^{+\infty}
|\hat{u}(x,t)|^2 dt
\end{equation*}
have the same value, finite or infinite.

Since $u(x,\cdot)$  belongs to $L^2(\RR)$ for almost every
$x\in\RR$, the same holds for $\hat{u}(x,\cdot)$ and, furthermore,
we have
\begin{equation*}
\ints_{-\infty}^{+\infty} |\hat{u}(x,t)|^2 dt = 2\pi
\ints_{-\infty}^{+\infty} |u(x,z)|^2 dz\quad \textmd{a.e. } x\in\RR.
\end{equation*}
By integrating the above equation and using Fubini-Tonelli's
theorem, we obtain
\begin{equation*}
\begin{split}
\ints_{\RR^2} |\hat{u}(x,t)|^2 dx dt & = \ints_{-\infty}^{+\infty}
dx \ints_{-\infty}^{+\infty} |\hat{u}(x,t)|^2 dt = 2\pi
\ints_{-\infty}^{+\infty} dx \ints_{-\infty}^{+\infty} |u(x,z)|^2 dz
\\ & =2 \pi \ints_{\RR^2} |u(x,z)|^2 dx dz < + \infty.
\end{split}
\end{equation*}
Therefore $\hat{u}(\cdot,t) \in L^2(\RR)$ for almost every
$t\in\RR$.

From \eqref{helm fourier in z}, it follows that
\begin{equation*}
\hat{u}(x,t) = a(t) \cos \sqrt{\lam-d^2} (x-h) + b(t) \sin
\sqrt{\lam-d^2} (x-h), \quad \textmd{for } x>h,
\end{equation*}
where $\lam=k^2n_*^2 -t^2$ and $d^2= k^2(n_*^2 - n_{cl}^2)$. Since
\begin{equation*}
\ints_{-\infty}^{+\infty} |\hat{u} (x,t)|^2 dx \geq
\ints_{h}^{+\infty} |\hat{u} (x,t)|^2 dx,
\end{equation*}
we obtain that $\hat{u} (x,t)$ can be {\it not identically zero}
only for some values $t=\lam_m^j \in (0, d^2]$ and, furthermore, in
that case
\begin{equation*}
\hat{u} (x,t) = a(t) v_s(x,\lam_s^m) + b(t) v_a(x,\lam_a^m).
\end{equation*}
Hence, for some $A, B\in\RR$ we should have
\begin{equation*}
u(x,z) = A Z_s(z) v_s(x,\lam_s^m) + B Z_a(z) v_a(x,\lam_a^m),
\end{equation*}
where $Z_j(z)=e^{\pm z \sqrt{k^2n_*^2 - \lam_m^j} },$ because $u$ is
a solution of \eqref{helm omogenea}. Since $u(x,\cdot) \in
L^2(\RR),$ then both $A$ and $B$ must be zero and hence $u\equiv 0$
on $\RR^2.$


\section{The spectrum-based solution satisfies the radiation
condition} \label{section 3}

It will be useful to introduce the following function
\begin{equation}\label{xh quadra}
[x]_h=
\begin{cases}
x+h, & x<-h, \\
0, & -h \leq x \leq h, \\
x-h, & x> h.
\end{cases}
\end{equation}

This section is devoted to the proof of the following result.

\vspace{1em}

\begin{figure}
\centering
\begin{minipage}[c]{0.4\textwidth}
\includegraphics[width=\textwidth]{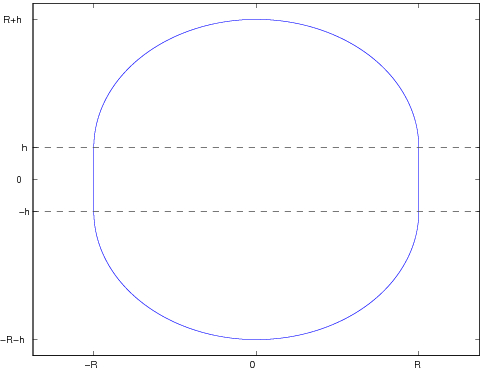}
\caption{The set $\Omega_R$.}  \label{fig omega0}
\end{minipage}
\hfill
\begin{minipage}{0.4\textwidth}
\includegraphics[width=\textwidth]{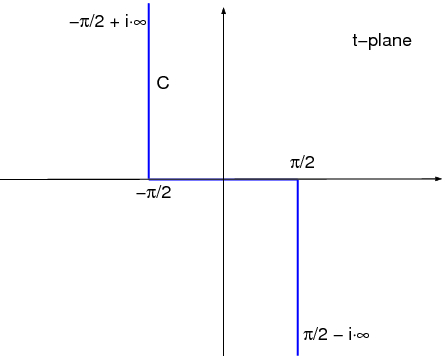}
\caption{The contour $\Cam$.} \label{Fig Cam}
\end{minipage}
\end{figure}

\begin{theorem}\label{teo verifica}
Let $f\in L^2(\RR^2)$ be such that $f\equiv 0 $ a.e. outside a
compact subset of $\RR^2$. Then, the spectrum-based solution
\eqref{u} of \eqref{helm} is the only solution in $L^{2,s}(\RR^2)$,
$s<-1$, such that
\begin{equation} \label{rad cond verifica}
\ints_0^\infty \ints_{\partial \Omega_R} \Big{|} \frac{\partial
u_0}{\partial \nu} - i\beta_0 u_0 \Big{|}^2 d\ell \, d R +
\sum_{l=1}^M \ints_0^\infty \ints_{\partial Q_R} \Big{|}
\frac{\partial u_l}{\partial \nu} - i\beta_l u_l \Big{|}^2 d\ell \,
d R < + \infty,
\end{equation}
where $\Omega_R$ is given by
\begin{equation} \label{omega rho}
  \Omega_R = \left\{ (x,z) \in \RR^2 : [x]_h^2 + z^2 \leq R^2  \right\}
\end{equation}
(see Fig. \ref{fig omega0}) and $Q_R= \left\{ (x,z) \in \RR^2 : |x|, |z| \leq R \right\}.$
\end{theorem}

\vspace{1em}

\begin{remark}\label{remark 2}
At the cost of extra computations, it may be proved
that Theorem \ref{teo verifica} also holds if we replace
\eqref{rad cond 2} by the more compact condition
\eqref{rad cond} with $\Omega_R$ given by \eqref{omega rho}.
\end{remark}

\vspace{1em}

We shall break the proof of Theorem \ref{teo verifica} up into three
steps. First, in \S \ref{subsection 3.1}, we will derive a handier
representation of the radiating part $G^{rad}$ of the Green's
function, as a suitable contour integral (see Lemma \ref{lemma
3.1}). Then, in \S \ref{subsection 3.2}, we will prove a uniform
asymptotic expansion for the quantity $\displaystyle\frac{\de G^{rad}}{\de \nu} -
i\beta_0 G^{rad}$ on the sets $\de \Omega_R$. Such an expansion will
be used in \S \ref{subsection 3.3} to carry out the proof of Theorem
\ref{teo verifica}, where we also test the radiation condition on
the guided components of $u$.


\subsection{Representing $G^{rad}$ as a contour integral}
\label{subsection 3.1}

We introduce the following functions:
\begin{equation*}
\{ x \}_h= x-[x]_h,
\end{equation*}
with $[x]_h$ given by \eqref{xh quadra}, and, for $\tau \in \CC$,
\begin{equation}\label{Phi_j}
\Phi_j(x,\tau) = \phi_j(\{x\}_h, d^2 + \tau^2) + \frac{
\phi_j'(\{x\}_h, d^2 + \tau^2)}{i\tau}, \quad \ja.
\end{equation}

With these notations, \eqref{vj ch2} and \eqref{sigma ch2} take the
more compact forms:
\begin{equation}\label{vj [x]h}
v_j(x,d^2+\tau^2)= \frac{1}{2} \Big\{ \Phi_j(x,\tau) e^{i\tau [x]_h}
+  \Phi_j(x,-\tau) e^{-i\tau [x]_h} \Big\}
\end{equation}
and
\begin{equation}\label{sigmaj [x]h}
\sigma_j(d^2+\tau^2) = \frac{1}{\Phi_j(h,\tau) \Phi_j(h,-\tau)},
\end{equation}
for $\ja$.

\vspace{1em}

\begin{lemma} \label{lemma 3.1}
Let $\Cam$ be the contour from $-\frac{\pi}{2} + i\cdot \infty$ to
$\frac{\pi}{2} - i\cdot \infty$ shown in Fig. \ref{Fig Cam} and let
$G^{rad}$ be the function in \eqref{G^rad ch2}. Then,
\begin{equation*}
G^{rad}=\somma \ints_{\Cam} \Big[ A_j^+(x,\xi;t) e^{i\beta_0
\alpha_+(x,z;\xi,\zi;t)} + A_j^-(x,\xi;t) e^{i\beta_0
\alpha_-(x,z;\xi,\zi;t)} \Big] dt,
\end{equation*}
with
\begin{equation*} 
A_j^\pm (x,\xi;t)= \frac{1}{8\pi i} \! \cdot \!
\frac{\Phi_j(x,\beta_0 \sin t) \Phi_j(\xi,\pm \beta_0 \sin
t)}{\Phi_j(h,\beta_0 \sin t) \Phi_j(h,-\beta_0 \sin t)},
\end{equation*}
and
\begin{equation*} 
\alpha_\pm (x,z;\xi,\zi;t) = ( [x]_h \pm [\xi]_h ) \sin t + |z-\zi|
\cos t,
\end{equation*}
$t\in\CC$, and where $\Phi_j$, $\ja$, is given by \eqref{Phi_j}. In
particular, the following equivalent expression for $G^{rad}$ will
also be useful:
\begin{equation}\label{Grad con la g}
G^{rad} = \ints_{\Cam} g(x,\xi;t) \, e^{i\beta_0 ([x]_h \sin t +
|z-\zi| \cos t)} dt,
\end{equation}
where
\begin{equation}\label{g xh xi t}
g(x,\xi;t)= \somma \left[ A_j^+ (x,\xi;t) e^{i[\xi]_h \sin t} + A_j^-
(x,\xi;t) e^{-i[\xi]_h \sin t}\right].
\end{equation}
(Notice that $g$ does not depend on $x$ for $|x|\geq h$.)
\end{lemma}

\begin{proof}
We first take \eqref{G^rad ch2} and make the change of variable
$\tau=\sqrt{\lam - d^2}$ to obtain:
\begin{equation*}
G^{rad} =  \frac{1}{4\pi i} \sum_{\ja} \ints_{-\infty}^{+\infty}
\frac{e^{i|z-\zi| \sqrt{\beta_0^2 - \tau^2}}}{\sqrt{\beta_0^2 -
\tau^2}}  v_j(x,\tau^2+d^2) v_j(\xi,\tau^2+d^2) \sigma_j(\tau^2+d^2)
d\tau;
\end{equation*}
here, we also used the fact that all the relevant quantities subject
to integration are even functions of $\tau$. With the help of
\eqref{vj [x]h} and \eqref{sigmaj [x]h}, and simple manipulations,
we can infer that
\begin{multline*}
G^{rad} =  \frac{1}{8\pi i} \sum_{\ja}\: \ints_{-\infty}^{+\infty}
\bigg\{ \frac{\Phi_j(x,\tau) \Phi_j(\xi,\tau) }{ \Phi_j(h,\tau)
\Phi_j(h,-\tau)} e^{i \big[\tau([x]_h+[\xi]_h) + |z-\zi|
\sqrt{\beta_0^2- \tau^2}\, \big]} \\ + \frac{\Phi_j(x,\tau)
\Phi_j(\xi,-\tau) }{ \Phi_j(h,\tau) \Phi_j(h,-\tau)} e^{i
\big[\tau([x]_h - [\xi]_h) + |z-\zi| \sqrt{\beta_0^2- \tau^2}\,
\big]} \bigg\} d\tau.
\end{multline*}
The conclusion is then readily obtained by splitting the interval of
integration up into the three intervals
$(-\infty,-\beta_0),[-\beta_0,\beta_0]$ and $(\beta_0,+\infty)$ and
by subsequently making the change of variable $\tau=\beta_0 \sin t$, with
$t\in \Cam$.
\end{proof}

\vspace{2em}

\begin{lemma} \label{lemma derivate}
For every $\xi,\zi$ fixed, we have:
\begin{subequations} \label{de Grad gradiente}
\begin{equation}\label{de Grad de x}
\frac{\de G^{rad}}{\de x} = i\beta_0 \ints_{\Cam} g(h\sign{x},\xi;t)
\sin{t} \, e^{i\beta_0 ([x]_h \sin t + |z-\zi| \cos t)} dt,
\end{equation}
for $|x|\geq h$ and $z\neq \zi$;
\begin{equation}\label{de Grad de z}
\frac{\de G^{rad}}{\de z} = i\beta_0 \sign{(z-\zi)} \ints_{\Cam}
g(x,\xi;t) \cos{t} \, e^{i\beta_0 ([x]_h \sin t + |z-\zi| \cos t)}
dt,
\end{equation}
\end{subequations}
for $z\neq \zi$.

In particular, on the set
$(0,\zi) + \de \Omega_R$ given by \eqref{omega rho}, we have:
\begin{subequations} \label{de Grad de nu}
\begin{equation}\label{de Grad de nu 1}
\frac{\de G^{rad}}{\de \nu} - i\beta_0 G^{rad} = i\beta_0
\ints_{\Cam} g(x,\xi;t) [\cos{t}-1] \, e^{i\beta_0
R\cos t} dt,
\end{equation}
for $z-\zi=R$ and $|x| \leq h$, and
\begin{equation}\label{de Grad de nu 2}
\frac{\de G^{rad}}{\de \nu} - i\beta_0 G^{rad} = i\beta_0
\ints_{\Cam} g(h,\xi;t)
[\cos{(t-\vt)}-1] \, e^{i\beta_0 R
\cos{(t-\vt)}} dt,
\end{equation}
\end{subequations}
where $\nu$ is the normal to $(0,\zi) + \de
\Omega_R$ and we have set $[x]_h= R \sin \vt$ and $z-\zi=R \cos \vt$ with
$\vt\in[0,\pi/2)$ and $R> h.$

Formulas analogous to \eqref{de Grad de nu} hold for the
remaining values of $\vt$ in $[-\pi, \pi).$
\end{lemma}

\vspace{0.5em}

\begin{proof}
Since $z\neq \zi$ and $\IM ([x]_h \sin t + |z-\zi| \cos{t}) \to
+\infty$ as $t\to \infty$ on $\Cam$, the integrands in \eqref{de
Grad de x} and \eqref{de Grad de z} vanish exponentially as $t\to
\infty$ on $\Cam$, since $g$ is bounded (see Lemma \ref{lemma g
bounded}). Thus, \eqref{de Grad de x} and \eqref{de Grad de z}
follow from an application of Lebesgue's dominated convergence
Theorem.
\end{proof}


\subsection{Uniform asymptotic estimates for $\displaystyle\frac{\de
G^{rad}}{\de \nu} - i\beta_0 G^{rad}$}
\label{subsection 3.2}

Aiming to estimate, as $R\to\infty$, the function $\frac{\de
G^{rad}}{\de \nu} - i\beta_0 G^{rad}$ given by \eqref{de Grad de
nu}, we need to deform the contour $\Cam$ to a more convenient one.

Without loss of generality we can assume that $\vt \in [0,\pi/2]$.
We define the new contour $\Cam_\vt$ (see Fig. \ref{Fig
Cam_deformato}) as follows:
\begin{equation*}
\Cam_\vt = \bigcup_{j=1}^5 \Gamma_j,
\end{equation*}
where
\begin{equation*} 
\delta_1=\arccos\frac{2\beta_0}{\sqrt{4\beta^2_0+d^2-\gamma_M}},
\quad \delta_2=\arcsinh \frac{\sqrt{d^2-\gamma_M}}{2\beta_0},
\end{equation*}
(notice that $\cos\delta_1 \cosh\delta_2=1$) and
\begin{eqnarray*}
&& \Gamma_1 = \{ t=t_1+it_2 \in\CC:\  \RE (\cos{t})=1,\ \IM(\cos{t})
\geq 0, \ -\frac{\pi}{2} < t_1 \leq
  -\delta_1,\ t_2\geq \delta_2 \} ,\\
&& \Gamma_2= \{ t\in \CC:\  -\delta_1 \leq t_1 \leq -\delta_1+\vt,
t_2=\delta_2\}, \\
&& \Gamma_3= \{ t\in \CC:\ \RE[\cos{(t-\vt)}]=1,\ \IM[
\cos{(t-\vt)}]\geq
  0,\ |t_1-\vt| \leq \delta_1,\ |t_2|\leq \delta_2 \} ,\\
&& \Gamma_4 = \{ t\in \CC:\  \delta_1+\vt \leq t_1 \leq \pi
-\delta_1 , \ t_2=-\delta_2 \}, \\
&& \Gamma_5 = \{ t\in \CC:\ \RE(\cos{t})=-1,\ \IM (\cos{t}) \geq 0,\
\frac{\pi}{2} < t_1 \leq \pi-\delta_1,\ t_2 \leq -\delta_2\}.
\end{eqnarray*}

This choice of $\Cam_\vt$ is suggested by the following three
remarks:
\begin{enumerate}[(i)]
\item $\Cam \cup \Cam_\vt$ does not contain in its interior the
poles of $g$ (which correspond to the guided part \eqref{G^g ch2} of
$G$);

\item $\Gamma_3$ is part of the steepest descent path of
$\cos(t-\vt)$;

\item $\Gamma_1,\Gamma_2,\Gamma_4,\Gamma_5$ are chosen to complete
the contour $\Cam \cup \Cam_\vt$ and to fulfill Lemma \ref{lemma
gamma3 piu O} below.

\end{enumerate}

\begin{figure}
\begin{center}
\includegraphics[width=0.8\textwidth]{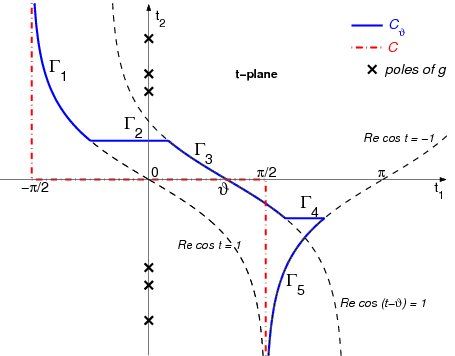}
\caption{The contour $\Cam$.}\label{Fig Cam_deformato}
\end{center}
\end{figure}

\vspace{1em}

By (i), it is clear that we can write
\begin{equation*}
\frac{\de G^{rad}}{\de \nu} - i\beta_0 G^{rad} = i\beta_0
\ints_{\Cam_0} g(x,\xi;t) [\cos{t}-1] \, e^{i\beta_0 R\cos t} dt,
\end{equation*}
for $|x| \leq h$, and
\begin{equation*}
\frac{\de G^{rad}}{\de \nu} - i\beta_0 G^{rad} = i\beta_0
\ints_{\Cam_\vt} g(h,\xi;t) [\cos{(t-\vt)}-1] \, e^{i\beta_0 R
\cos{(t-\vt)}} dt,
\end{equation*}
for $x \geq h$.

\vspace{1em}

\begin{lemma}\label{lemma gamma3 piu O}
Let $(x,z)\in (0,\zi)+ \de \Omega_R$. The following estimates hold
for $R\to \infty$:
\begin{subequations} \label{eq gamma3 piu O}
\begin{equation} \label{eq gamma3 piu O x<h}
\frac{\de G^{rad}}{\de \nu} - i\beta_0 G^{rad} = i\beta_0
\ints_{\Gamma_3} g(x,\xi;t) [\cos{t}-1] \, e^{i\beta_0 R \cos{t}} dt
+ \OO(e^{-c\beta_0 R}),
\end{equation}
for $|x|\leq h$, and
\begin{equation} \label{eq gamma3 piu O x>h}
\frac{\de G^{rad}}{\de \nu} - i\beta_0 G^{rad} = i\beta_0
\ints_{\Gamma_3} g(h,\xi;t) [\cos{(t-\vt)}-1] \, e^{i\beta_0 R
\cos{(t-\vt)}} dt  + \OO(e^{-c\beta_0 R}),
\end{equation}
\end{subequations}
for $x \geq h$, $\vt \in [0,\pi/2]$, where
\begin{equation*}
c= \sqrt{\frac{d^2- \gamma_M}{4 \beta_0^2 + d^2 - \gamma_M}}
\cdot\min \left( 1,\frac{\sqrt{d^2-\gamma_M}}{2\beta_0} \right).
\end{equation*}
\end{lemma}

\begin{proof}
We shall prove only \eqref{eq gamma3 piu O x>h} since \eqref{eq
gamma3 piu O x<h} follows analogously. We preliminarily observe that
\begin{equation}\label{exp -R}
\IM \cos(t-\vt) \geq c,
\end{equation}
for $t \in \Gamma_1,\Gamma_2,\Gamma_4,\Gamma_5$ and $\vt \in
[0,\pi/2]$. From \eqref{exp -R}, we easily obtain that
\begin{equation*}
\Bigg{|} \ints_{\Gamma_j} g(h,\xi;t) [\cos(t-\vt) - 1] e^{i\beta_0 R
\cos{(t-\vt)}} dt \Bigg{|} \leq \frac{K \pi}{2}
(\cosh{\delta_2}+1)e^{-c \beta_0 R }, \quad j=2,4,
\end{equation*}
where $K$ is a bound for $g$ (see Lemma \ref{lemma g bounded}).
Thus, it remains to prove that
\begin{equation*}
\ints_{\Gamma_j} g(h,\xi;t) [\cos{(t-\vt)}-1] \, e^{i\beta_0 R
\cos{(t-\vt)}} dt = \OO(e^{-c\beta_0 R}), \quad j=1,5,
\end{equation*}
uniformly as $R\to \infty$, for $\vt \in [0,\pi/2]$. We carry out
the details for $j=1$, the case $j=5$ is completely analogous. We
first use Lemma \ref{lemma g asint} to write that
\begin{equation*}
\ints_{\Gamma_1} g(h,\xi;t) [\cos{(t-\vt)}-1] \, e^{i\beta_0 R
\cos{(t-\vt)}} dt = J(R) + \OO(e^{-c\beta_0 R}),
\end{equation*}
since \eqref{exp -R} holds; here,
\begin{equation*}
J(R) = i\beta_0 \ints_{\Gamma_1} \bigg( 1+ \frac{i}{2\beta_0 \sin t}
\ints_{\{\xi\}_h}^h p(y) dy \bigg) [\cos{(t-\vt)}-1] e^{i\beta_0 [
R\cos(t- \vt) + (h-\xi) \sin t]} dt,
\end{equation*}
with $p(y)=d^2 - q(y)$. Let $\psi(t)= R\cos(t- \vt) + (h-\xi) \sin
t$ and $\delta=\delta_1+i\delta_2$; an integration by parts yields
\begin{multline*}
J (R) =\frac{e^{i\beta_0 \psi(\delta )}}{\psi'(\delta)}
\bigg(1+\frac{i}{2\beta_0 \sin \delta} \ints_{\{\xi\}_h}^h
p(y) dy \bigg)[1- \cos{(\delta -\vt)}] \\
+ \ints_{\Gamma_1} \frac{e^{i\beta_0 \psi(t)}}{\psi'(t)^2} \bigg\{
[\cos{(t-\vt)}-1] \bigg[R- \frac{i \psi(t)}{2\beta_0 \sin t}
\ints_{\{\xi\}_h}^h p(y) dy \bigg] \\ + (h-\xi)(\sin{t} - \sin{\vt})
- \frac{i \psi'(t) (\cos{t} - \cos{\vt})}{2\beta_0 \sin^2 t}
\ints_{\{\xi\}_h}^h p(y) dy \Bigg\} dt.
\end{multline*}
From \eqref{exp -R} and since
\begin{eqnarray*}
& \sinh t_2 \leq |\cos t|,|\sin t| \leq \cosh t_2,\\
&|\psi(t)| \leq \beta_0 ( R + |h-\xi| ) \cosh{t_2},\quad
|\psi'(t)| \geq \frac{1}{2} \beta_0 R \sinh{t_2} ,\quad
\textmd{for } R \geq 2 |h-\xi| \coth{\delta_2} ,
\end{eqnarray*}
for $t\in \Gamma_1$, we obtain that $J(R) = \OO(e^{-c\beta_0 R})$,
as $R\to \infty$.
\end{proof}

\vspace{2em}

\begin{theorem}\label{teo Grad stima uniforme}
On $\de \Omega_R$, we have
\begin{equation} \label{Grad stima cond rad}
\frac{\de G^{rad}}{\de \nu} - i \beta_0 G^{rad} = \OO \left(
R^{-\frac32} \right),
\end{equation}
uniformly as $R\to \infty$.
\end{theorem}

\begin{proof}
First, we estimate the left-hand side of \eqref{Grad stima cond rad}
on the sets $(0,\zi)+\de \Omega_R$. By Lemma \ref{lemma gamma3 piu
O}, we only need to estimate the first addendum in \eqref{eq gamma3
piu O}. We prove the estimate for \eqref{eq gamma3 piu O x>h}; the
estimate for \eqref{eq gamma3 piu O x<h} follows exactly in the same
way.

Since $\Gamma_3$ is part of the steepest descent path, the steepest
descent method (see \cite{BO}) suggests to change the variables in
the first addendum in \eqref{eq gamma3 piu O x>h}: by setting
$\cos(t-\vt)=1+iy^2$, we obtain
\begin{equation*}
\ints_{\Gamma_3} g(h,\xi;t) [\cos{(t-\vt)}-1] \, e^{i\beta_0 R
\cos{(t-\vt)}} dt = - 4 i e^{i\beta_0 R} \ints_{0}^{y_0} y^2
e^{-\beta_0 R y^2} \frac{g(h,\xi;t(y))}{\sqrt{y^2 - 2i }} dy,
\end{equation*}
with $y_0=(\sin\delta_1 \sinh\delta_2)^{\frac{1}{2}}$. Thanks to
Lemma \ref{lemma g bounded},
\begin{equation*}
\bigg{|}\ints_{\Gamma_3} g(h,\xi;t) [\cos{(t-\vt)}-1] \, e^{i\beta_0
R \cos{(t-\vt)}} dt \bigg{|} \leq 2\sqrt{2} K \ints_0^{y_0} y^2
e^{-\beta_0 R y^2} dy \leq K \sqrt{\frac{\pi}{2 \beta_0^3}}\
R^{-\frac32},
\end{equation*}
where $K$ is a bound of $g$. Therefore, \eqref{eq gamma3 piu O x>h}
implies that
\begin{equation*}
\bigg{|} \frac{\de G^{rad}}{\de \nu} - i \beta_0 G^{rad} \bigg{|}
\leq K \sqrt{\frac{\pi}{2\beta_0 }}\ R^{-\frac32},
\end{equation*}
on the sets $(0,\zi)+\de \Omega_R$.

By using exactly the same argument as before, we can prove that the
derivatives of $G^{rad}$ are $\OO (R^{-\frac{1}{2}})$ on the sets
$(0,\zi)+\de \Omega_R$, uniformly as $R \to \infty$; we reach the
conclusion \eqref{Grad stima cond rad} by observing that $\nu_{\de
\Omega_R} - \nu_{(0,\zi)+ \de \Omega_R} = \OO (R^{-1})$, as $R \to
\infty$.
\end{proof}

\vspace{2em}

\subsection{Proof of Theorem \ref{teo verifica}} \label{subsection 3.3}
Since $f\in L^2(\RR^2)$ and $f$ has compact support, from Corollary
5.1 in \cite{CM} we have that $u\in L^{2,s}(\RR^2)$, $s<-1$.

Thus, it remains to prove that \eqref{u} satisfies \eqref{rad cond}.
In order to do it, we shall check the following facts:
\begin{enumerate}[(i)]
\item if $u$ is given by \eqref{u} and $u_l$, $l=1,\ldots,M$, is
computed via \eqref{u_l}, the remainder part $u_0$ of
$u$, given by \eqref{u_0}, equals the function $u^{rad}$ in \eqref{u
rad originale};

\item  $u$ satisfies \eqref{rad cond verifica}.
\end{enumerate}

We preliminarily notice that
\begin{equation*} \ints_0^{+\infty}
\ints_{\de \Omega_R} \bigg{|} \frac{\de u_0}{\de \nu} - i \beta_0
u_0 \bigg{|}^2 d\ell dR < + \infty
\end{equation*}
is easily verified thanks to Theorem \ref{teo Grad stima uniforme}.

The following property of orthogonality is useful to check (i).

\vspace{1em}

\begin{lemma} \label{lemma ortogonali}
Let $e(x,\gamma_l)$, $l=1,\ldots, M$, and $v_j(x,\lam)$, $\ja$, be the
solutions of \eqref{eq v ch2} given by \eqref{e autof} and \eqref{vj
ch2}, respectively, with $\lam > 0$. If $\lam \neq \gamma_l$, then
\begin{equation*}
    \ints_{-\infty}^{+\infty} e(x,\gamma_l) v_j(x,\lam) dx = 0.
\end{equation*}
\end{lemma}
\begin{proof}
We multiply the following equations
\begin{eqnarray*}
&& e''(x,\gamma_l) + [\gamma_l - q(x)] e(x,\gamma_l) = 0, \\
&& v''(x,\lam) + [\lam - q(x)] v(x,\lam) = 0,
\end{eqnarray*}
by $v(x,\lam)$ and $e(x,\gamma_l)$, respectively, and integrate in
$x$ over an interval $(a,b)$. An integration by parts gives:
\begin{equation*}
\begin{split}
(\gamma_l-\lam) \ints_a^b  e(x,\gamma_l) v(x,\lam) dx & = \ints_a^b
\left[  e(x,\gamma_l) v''(x,\lam) - e''(x,\gamma_l)
v(x,\lam) \right] dx \\
& = \left[  e(x,\gamma_l) v'(x,\lam) -  e'(x,\gamma_l) v(x,\lam)
\right]_a^b.
\end{split}
\end{equation*}
The conclusion follows by observing that $ e(x,\gamma_l) $ and its
first derivative vanish exponentially as $|x|\to \infty$, while
$v(x,\lam)$ and $v'(x,\lam)$ are bounded.
\end{proof}

\vspace{1em}

Now, by \eqref{u}, \eqref{u_l} and Lemma \ref{lemma ortogonali}, we
have that
\begin{equation*}
u_l(x,z)= \ints_{\RR^2} G_l^g (x,z;\xi,\zi) f(\xi,\zi) d\xi d\zi,
\quad l=1,\ldots,M,
\end{equation*}
with $G_l^g$ given by \eqref{G l guid} and thus $u_0=u^{rad}$.

To complete the proof it remains to check (ii) for $l=1,\ldots,M$.
When  $z$ is large enough, we have
\begin{equation*}
\frac{\de u_l}{\de \nu} - i \beta_l u_l = 0,\quad l=1,\ldots,M,
\end{equation*}
on $\de Q_R \cap \{(x,z) : |z|=R \}$, since $\frac{\de}{\de \nu} =
\pm \frac{\de}{\de z}$. Thanks to \eqref{vj guidate}, we easily find
that
\begin{equation*}
\bigg{|} \frac{\de u_l}{\de \nu} - i \beta_l u_l \bigg{|} = \OO
\left( e^{-R \sqrt{d^2-\gamma_l}} \right),
\end{equation*}
as $R\to \infty$ on $\de Q_R \cap \{(x,z) : |x|=R \}$ and thus we
obtain that
\begin{equation*} \ints_0^{+\infty}
\ints_{\de Q_R} \bigg{|} \frac{\de u_l}{\de \nu} - i \beta_l u_l
\bigg{|}^2 d\ell dR < + \infty, \quad l=1,\ldots,M,
\end{equation*}
which completes the proof.

\appendix

\section{Asymptotic Lemmas}
In what follows, $BV(\RR)$ denotes the space of functions with
bounded variation.

\begin{lemma} \label{lemma Phi asint}
Let $T$ be a non-negative number, $q\in BV(\RR)$ and
\begin{equation*} 
p(x)=d^2 - q(x),\quad x\in \RR.
\end{equation*}
Then, the following asymptotic estimates for the functions $\Phi_s$
and $\Phi_a$ given by \eqref{Phi_j} hold uniformly as $|\tau| \to +\infty$, for $x\in \RR$ and $|\IM
\tau | \leq T$:
\begin{eqnarray}
&& \Phi_s(x,\tau) = \Bigg[ 1 + \frac{i}{2\tau} \ints_0^{\{x\}_h}
p(y) dy
\Bigg] e^{i\tau \{x\}_h} + \OO \bigg(\frac{1}{|\tau|^2} \bigg),  \label{Phi s asint}\\
&& \Phi_a(x,\tau) = \frac{\sqrt{\tau^2 + d^2}}{i \tau}\ \Bigg[ 1 +
\frac{i}{2\tau} \ints_0^{\{x\}_h} p(y) dy \Bigg] e^{i\tau \{x\}_h} +
\OO \bigg(\frac{1}{|\tau|^2} \bigg). \label{Phi a asint}
\end{eqnarray}

\end{lemma}

\begin{proof}
(i) First, we prove an estimate for $\phi_j$, $\ja$. From \eqref{eq
v ch2}, we know that $\phi_j$ satisfies
\begin{equation*}
\phi_j''(y,\lam) + [\tau^2 + p(y)] \phi_j(y,\lam) = 0,\quad y\in[-h,h].
\end{equation*}
We multiply the above equation by $\sin\tau(x-y)$,
integrate by parts twice and obtain the following integral equation:
\begin{equation} \label{eq integrale phi}
\phi_j(x,\lam) = \frac{\phi_j'(0,\lam)}{\tau} \sin{\tau x} +
\phi_j(0,\lam) \cos{\tau x} - \frac{1}{\tau} \ints_0^x p(y)
\sin(\tau(x-y)) \phi_j(y,\lam)\: dy.
\end{equation}
We set $\eta_j(x,\lam)= \sup\limits_{s\in [0,x]} |\phi_j(s,\lam)|$.
Since $|\sin \tau x|, |\cos \tau x| \leq \cosh{\tau_2 x}$
($\tau_2=\IM \tau$), from the above equation we have that
\begin{equation*}
\eta_j(x,\lam) \leq \left[ \frac{|\phi_j'(0,\lam)|}{|\tau|} +
|\phi_j(0,\lam)| \right] \cosh{\tau_2 x} + \frac{1}{|\tau|}
\ints_0^x p(y) \cosh\tau_2(x-y)\ \eta_j(y,\lam)\: dy,
\end{equation*}
and, by Gronwall's Lemma (see \cite{SC}), we get
\begin{multline*}
\eta_j(x,\lam) \leq  \left[ \frac{|\phi_j'(0,\lam)|}{|\tau|} +
|\phi_j(0,\lam)| \right] e^{\frac{1}{|\tau|} \ints_0^x p(y)
\cosh\tau_2(x-y)\: dy} \times \\ \times \left\{ 1 + \tau_2 \ints_0^x
e^{-\ints_0^s p(y) \cosh\tau_2(x-y)\: dy} \sinh{\tau_2 s}\ ds
\right\}.
\end{multline*}
Since $0\leq p(y) \leq d^2$, we have that
\begin{equation*}
\eta_j(x,\lam) \leq  \left[ \frac{|\phi_j'(0,\lam)|}{|\tau|} +
|\phi_j(0,\lam)| \right] \cosh{\tau_2 x} \: \exp\left\{\frac{d^2
\sinh{\tau_2 x}}{|\tau| \tau_2}\right\} .
\end{equation*}
If we assume $|\tau| \geq d$ and $x\in[-h,h]$, we finally get
\begin{equation}\label{phi stima}
|\phi_j(x,\lam)| \leq \sqrt{2} \cosh\tau_2 h\: \exp \left\{\frac{d
\sinh{T h}}{T} \right\}, \quad \ja.
\end{equation}

(ii) Now we prove \eqref{Phi s asint} and \eqref{Phi a asint}. Let
$q\in C^1(\RR)$. From \eqref{Phi_j}, by straightforward
manipulations we have:
\begin{equation} \label{Phi eq diff}
\Phi_j'(x,\tau) - i \tau \Phi_j(x,\tau) = \frac{i}{\tau}\: p(x)\:
\phi_j(x,\lam);
\end{equation}
by multiplying the above equation by $e^{-i\tau x}$, integrating by
parts twice and observing that
\begin{multline*}
2 \ints_0^x e^{-i\tau y} p(y) \phi_j(y,\lam) dy = \ints_0^x
e^{-i\tau y} p(y) \Phi_j(y,\lam) dy \\+ \frac{1}{i\tau} \ints_0^x
e^{-i\tau y} p'(y) \phi_j(y,\lam) dy - \left[ \frac{e^{-i\tau
y}}{i\tau} p(y) \phi_j(y,\lam) \right]_{y=0}^{y=x},
\end{multline*}
it follows that $\Phi_j$ satisfies
\begin{multline} \label{eq2bis lemma Phi}
\Phi_j(x,\tau) e^{-i\tau x} = \Phi_j(0,\tau) + \frac{i}{2\tau}
\ints_0^x e^{-i\tau y} p(y) \Phi_j(y,\lam) dy \\+ \frac{1}{2\tau^2}
\left\{ \ints_0^x e^{-i\tau y} p'(y) \phi_j(y,\lam) dy - e^{-i\tau
x} p(x) \phi_j(x,\lam) + p(0) \phi_j(0,\lam) \right\}.
\end{multline}

By setting $M_j(x,\tau) = \sup\limits_{s\in[0,x]} |\Phi_j(s,\lam)
e^{-i\tau s}| $ and from \eqref{phi stima}, we get
\begin{multline*}
M_j(x,\tau) \leq |\Phi_j(0,\tau)| + \frac{1}{2|\tau|} \ints_0^x p(y)
M_j(y,\tau) dy
\\ + \frac{1}{2|\tau|^2} \left\{ C \ints_0^x e^{\tau_2 y} |p'(y)| dy + C d^2 e^{\tau_2 x}
+ p(0)| \phi_j(0,\lam)| \right\},
\end{multline*}
for $|\tau| \geq d$, where $C$ is the right-hand side of \eqref{phi
stima}. Thus, Gronwall's Lemma yields the following estimate for
$M_j$:
\begin{multline*}
M_j(x,\tau) \leq \left[ |\Phi_j(0,\tau)| + \frac{C d^2}{2|\tau|^2} +
p(0) |\phi_j(0,\lam)| \right] \exp\left\{ \frac{1}{2|\tau|}
\ints_0^x p(y) dy \right\} \\ + \frac{C}{2|\tau|^2} \ints_0^x
\exp\bigg\{ \frac{1}{2|\tau|} \ints_s^x p(y) dy \bigg\} e^{\tau_2 s}
[\: |p'(s)| + \tau_2 d^2] ds ;
\end{multline*}
since $\Phi_s(0,\tau)=1$,
$\Phi_a(0,\tau)=\frac{\sqrt{\tau^2+d^2}}{i\tau}$ and $0\leq p(x)
\leq d^2$, we have
\begin{equation} \label{Mj stima}
M_j(x,\tau) \leq e^{\frac{dh}{2}} \left\{ 1+\sqrt{2} +
\frac{C}{2d^2} e^{T x} \left[ d^2 + |q|_{BV} \right] \right\},
\end{equation}
for $|\tau| \geq d$. By a standard approximation argument we can
infer that \eqref{Mj stima} holds for every $q \in BV(\RR)$. By
\eqref{eq2bis lemma Phi}, \eqref{phi stima} and \eqref{Mj stima}, we
get that
\begin{equation*}
\Phi_j(x,\tau) = \Phi_j(0,\tau) e^{i\tau x} + \OO
\Big(\frac{1}{|\tau|} \Big).
\end{equation*}
Again, from \eqref{eq2bis lemma Phi} and the above asymptotic
formula, we obtain \eqref{Phi s asint} and \eqref{Phi a asint}.
\end{proof}

\vspace{2em}

In Lemmas \ref{lemma g asint} and \ref{lemma g bounded}, we will use
the following inequality:
\begin{equation} \label{IM sin t limitata}
|\IM \sin t| \leq \max \bigg\{1, \frac{\sqrt{d^2-\gamma_M}}{2
\beta_0} \bigg\}, \quad t\in \Cam_{\vt},\ \vt\in [0,\pi/2].
\end{equation}

\vspace{2em}

\begin{lemma} \label{lemma g asint}
Let $g$ be defined by \eqref{g xh xi t}. Then, the following
asymptotic expansion
\begin{equation}\label{g asint}
g(x,\xi;t)= \frac{1}{4\pi i} \: e^{i \beta_0 ( \{x\}_h - \xi ) \sin
t} \bigg[ 1+ \frac{i}{2\beta_0 \sin t} \ints_{\{\xi\}_h}^{\{x\}_h}
p(y) dy \bigg] + \OO \bigg(\frac{1}{|\sin t|^2} \bigg)
\end{equation}
holds uniformly as $t \to \infty$ on $\Cam_\vt$ for $\vt\in
[0,\pi/2]$, $x \in \RR$ and $\xi$ bounded.
\end{lemma}

\begin{proof}
The proof is a straightforward consequence of Lemma \ref{lemma Phi
asint} and \eqref{IM sin t limitata}, and hence is omitted.
\end{proof}

\vspace{2em}

\begin{lemma} \label{lemma g bounded}
Let $g$ be given by \eqref{g xh xi t}. Then $g$ is a bounded
function of $x,\xi \in\RR$, if $\xi$ is bounded, and $t\in
\Cam_\vt$, $\vt \in [0,\pi/2]$.
\end{lemma}

\begin{proof}
(i) First, we prove an estimate for $\phi_j(x,\tau^2 + d^2)$ for
$|\tau| \leq d$, $|\IM \tau| \leq T$ and $|x| \leq h$. By setting
$\lam=\tau^2 + d^2$ and $\eta_j(x,\lam)= \sup\limits_{s\in [0,x]}
|\phi_j(s,\lam)|$ as before, from \eqref{eq integrale phi} and since $ \big{|}
\frac{\sin\tau x}{\tau x} \big{|}$ is bounded by the constant  $B=\sqrt{\cosh^2 (Th) +
\frac{\sinh^2(Th)}{(Th)^2}}$, we have
\begin{eqnarray*}
&& \eta_s(x,\lam) \leq \cosh(T x) + B\: \ints_0^{|x|} p(y)
|x-y|
\eta_s(y,\lam) dy, \\
&& \eta_a(x,\lam) \leq B\: \bigg\{ \sqrt{2} d
|x| + \ints_0^{|x|} p(y) |x-y| \eta_a(y,\lam) dy \bigg\};
\end{eqnarray*}
Gronwall's Lemma yields
\begin{equation}\label{phi bound <d}
|\phi_j(x,\lam)| \leq \min \left\{ \cosh(T h), \sqrt{2} dh B
\right\} \exp \left( B\: \frac{d^2h^2}{2}
\right),
\end{equation}
for $|x| \leq h$, $|\tau| \leq d$ and $|\IM \tau| \leq T$.

(ii) To complete the proof, we notice that from \eqref{Phi eq diff}
it follows that
\begin{equation*}
\bigg{|} \frac{\Phi_j(x,\tau)}{\Phi_j(h,\tau)} \bigg{|} \leq \frac{d
e^{Th}}{|\tau \Phi_j(h,\tau)|} \bigg( \sqrt{2} + dh e^{Th}
\sup\limits_{x\in [-h,h]} |\phi_j(x,\lam)| \bigg),
\end{equation*}
and since $\tau \Phi_j(h,\tau) \neq 0 $ far from the poles of $g$,
we have that $\frac{\Phi_j(x,\tau)}{\Phi_j(h,\tau)}$ is bounded for
$\tau=\beta_0 \sin t$, $t\in \Cam_{\vt}$, $\vt \in [0,\pi/2]$ and
for $|\tau|\leq d$. Thus, the assertion of the lemma follows from
\eqref{phi bound <d} and Lemma \ref{lemma g asint}.
\end{proof}

\end{document}